\documentclass[12pt,centertags]{amsart}
\usepackage[colorlinks,linkcolor=blue,citecolor=blue,bookmarksnumbered]{hyperref}
\usepackage{amsmath,amstext,amsthm,a4,amssymb,amscd}
\usepackage[mathscr]{eucal}
\usepackage{mathrsfs}
\usepackage{epsf}
\numberwithin{equation}{section}

\newtheorem{thm}{Theorem}[section]

\theoremstyle{definition}

\newtheorem{con}[thm]{Conjecture}
\theoremstyle{definition}

\theoremstyle{definition}
\newtheorem{defn}[thm]{Definition}
\newcommand{\be}{\begin{eqnarray}}
\newcommand{\ee}{\end{eqnarray}}

\newcommand{\comment}[1]{}

\begin{document}

\title[On the generalized Geroch conjecture for complete spin manifolds]{On the generalized Geroch conjecture for\\ complete spin manifolds}
\author{Xiangsheng Wang\ \ and\ \ Weiping Zhang}

\address{School of Mathematics, Shandong
University, Jinan, Shandong 250100, P.R. China}
\email{xiangsheng@sdu.edu.cn}

\address{Chern Institute of Mathematics \& LPMC, Nankai
University, Tianjin 300071, P.R. China}
\email{weiping@nankai.edu.cn}

\begin{abstract}  Let $W$ be a closed area  enlargeable manifold in the sense of    Gromov-Lawson and $M$ be a   noncompact spin manifold,  we show that the connected sum $M\# W$   admits no complete  metric of positive scalar curvature. When $W=T^n$, this provides a positive answer to the generalized Geroch conjecture in the spin setting. 
\end{abstract}

\maketitle

\setcounter{section}{-1}

\section{Introduction} \label{s0}

It has been an important subject in differential  geometry  to study when a smooth manifold carries a Riemannian metric of positive scalar curvature.  A famous theorem  of Gromov and Lawson \cite{GL80}, \cite{GL83} states that an area  enlargeable   manifold (in the sense of \cite{GL83})  does not carry a metric of positive scalar curvature.

\begin{defn}\label{t0.1} (Gromov-Lawson  \cite{GL83})
One calls a closed manifold $W$ (carrying a metric $g^{TW}$)  an area enlargeable manifold if for any $\epsilon>0$, there is a   covering manifold $\pi:\allowbreak\widehat W_\epsilon\rightarrow W$ (carrying the lifted metric), with $\widehat W_\epsilon$ being spin,  and a smooth map $f:\widehat W_\epsilon\rightarrow S^{\dim W}(1)$ (the standard unit sphere), which is constant near infinity and has non-zero degree, such that for any two form $\alpha\in \Omega^2(S^{\dim W}(1))$,  one has $|f^*(\alpha)|\leq \epsilon |\alpha|$.
\end{defn}

It is clear that the area enlargeability does not depend on the metric $g^{TW}$.

\begin{thm}\label{t0.2}  Let $W$ be a  closed area enlargeable manifold and $M$ an arbitrary  spin manifold of equal dimension, then the connected sum $M\# W$ does not admit any complete metric of positive scalar curvature.
\end{thm}

When   $M$ is closed, Theorem \ref{t0.2} is exactly the   Gromov-Lawson theorem \cite{GL80}, \cite{GL83}  mentioned at the begining. 
When $W=T^n$, Theorem \ref{t0.2} solves the following generalized Geroch conjecture (cf. \cite[Conjecture 1.4]{Zhu22}) in the spin setting.

\begin{con}\label{t0.3}  For any manifold $M$ of dimension $n$, there is no complete metric of positive scalar curvature on $T^n\# M$.
\end{con}

The remarkable fact (cf. \cite{Lo99}) is that     Conjecture \ref{t0.3}  for   the case of compact $M$ (first proved by Schoen-Yau \cite{SY79} in dimension $\leq 7$)  implies the     positive mass theorem for $M$, which in the spin case was proved by Witten \cite{W81} (see also \cite{PT})  using Dirac operators (the classical positive mass theorem in dimension three was first proved by Schoen-Yau \cite{SY79a}   using minimal hypersurface method  which works  for dimension $\leq 7$, cf. \cite{SY79b}), while a proof in the nonspin case for all dimensions  is given by Schoen-Yau \cite{SY17} by further developing  their   minimal hypersurface techniques (see also Lohkamp \cite{Lo16} for another minimal hypersurfaces approach in the higher dimensional situation).

 If $3\leq \dim M\leq 7$, then Conjecture \ref{t0.3} has been proved   for arbitrary $M$ by Chodosh and Li \cite{CL20}, with the case of $\dim M=3$ also proved by Lesourd-Unger-Yau \cite{LUY20}. A recent paper by Zhu \cite{Zhu22} shows that Conjecture \ref{t0.3} implies the positive mass theorem with arbitrary ends, which in the spin setting has been proved in \cite{BC03} (see also \cite[Theorem B]{CZ21}).  Thus Theorem \ref{t0.2} gives an alternate proof of the positive mass theorem with arbitrary ends in the spin setting.

  Our proof of Theorem \ref{t0.2} is based on  deformed Dirac operators as was used in \cite{Z20}. Indeed, by using a trick in \cite{SZ18} (which goes back to \cite{GL83}), we show that Theorem \ref{t0.2} reduces to the situation already considered in \cite{Z20}.

\section{A proof of Theorem \ref{t0.2}} \label{s1}

Let $W$ be a  closed area enlargeable manifold. Let $M$ be a   ${\rm spin}$   manifold. Without loss of generality, we assume that $W$ is spin    and that $\dim M=\dim W=n$.  Let   $h^{TW}$  be a metric on   $TW$.

As in \cite{SZ18}  (which goes back to \cite{GL83}), we fix a point $p\in W$. For any $r\geq 0$, let $B^W_p(r)=\{y\in W\,:\, d(p,y)\leq r\}$. Let $b_0>0$ be a fixed sufficiently small number. Then the connected sum $M\# W$ can be constructed so that  the hypersurface $\partial B^W_p(b_0)$, which is the boundary of $B^W_p(b_0)$,  cuts $M\# W$ into two parts: the part $W\setminus B^W_p(b_0)$ and the rest part coming from $M$ (by attaching the boundary of a ball in $M$ to $\partial B^W_p(b_0)$).

 For any $\epsilon>0$, let $\pi:\widehat W_\epsilon\rightarrow W$ be a  covering manifold verifying Definition \ref{t0.1}, carrying lifted geometric data from that of $W$. Let $b_0>0$ be small enough so that for any $p',\,q'\in\pi^{-1}(p)$ with $p'\neq q'$, $\overline{B^{\widehat W_\epsilon}_{p'}(4b_0)}\cap \overline{B^{\widehat W_\epsilon}_{q'}(4b_0)}=\emptyset$. It is clear that one can choose $b_0>0$ not depending on $\epsilon$.

Let $h: W\rightarrow W$ be a smooth map such that $h={\rm Id}$ on $W\setminus B^W_p(3b_0)$, while $h( B_p^W(2b_0))=\{p\}$.  It lifts to a map $\widehat h_\epsilon:\widehat W_\epsilon\rightarrow \widehat W_\epsilon$ verifying  that $\widehat h_\epsilon={\rm Id}$ on $\widehat W_\epsilon\setminus \bigcup_{p'\in\pi^{-1}(p)}B_{p'}^{\widehat W_\epsilon}(3b_0)$, while for any $p'\in\pi^{-1}(p)$, $\widehat h_\epsilon(b_{p'}^{\widehat W_\epsilon}(2b_0))=\{p'\}$.

Let $f :  \widehat W_\epsilon\rightarrow S^{n}(1)$ be the map given in Definition \ref{t0.1}, where for simplicity we assume that each $ \widehat W_\epsilon$ is compact. Set $\widehat f=f\circ \widehat h_\epsilon:\widehat{  W}_\epsilon\rightarrow S^n(1)$. Then ${\rm deg}(\widehat f)={\rm deg}(f)$ and  there is a positive constant $c>0$ (not depending on $\epsilon$) verifying that for any $\alpha\in\Omega^2(S^n(1))$, one has 
\begin{align}\label{1.1}
|\widehat f^*(\alpha)|\leq c\,\epsilon\, |\alpha|.
\end{align} 

The connected sum $M\# W$ lifts naturally to $\widehat W_\epsilon$ where near each $p'\in\pi^{-1}(p)$, we do the lifted connected sum, i.e., do the connected sum $M_{p'}\# B_{p'}^{\widehat W_\epsilon}(2b_0)$, where $M_{p'}$ is a copy of $M$. We denote the resulting manifold by $\widehat M\#\widehat W_\epsilon$. 
Cleary, any metric on $T(M\# W)$ lifts to a metric on $T(\widehat M\#\widehat W_\epsilon)$. 

We extend $\widehat f:\widehat W_\epsilon\rightarrow S^n(1)$ to $ \widehat M\#\widehat W_\epsilon$ by setting $\widehat f(M_{p'}\# B_{p'}(2b_0))=f(p')$ for any $p'\in\pi^{-1}(p)$. Clearly,  $\widehat f:\widehat M\#\widehat W_\epsilon\rightarrow S^n(1)$ is locally constant on $(\widehat M\#\widehat W_\epsilon)\setminus  (\widehat W_\epsilon \setminus \cup _{p'\in\pi^{-1}(p)} B_{p'}(b_0))$. Moreover, we still have ${\rm deg}(\widehat f)={\rm deg}(f)$.

Let $g^{T(M\# W)}$ be any complete metric on $T(M\# W)$ of positive scalar curvature. 
Since $ \overline{ W  \setminus   B_{p}(b_0)}$ is compact in 
$  M\#  W  $,   there is $\delta>0$ such that the corresponding scalar curvature satisfies that
\begin{align}\label{1.2}
 k^{g^{T(M\# W)}}\geq \delta\ \ {\rm on}\ \   (W\setminus B_{b_0}(p))\subset M\# W.
\end{align} 
A similar inequality  also holds for the scalar curvature of the  lifted metric  $g^{T(\widehat M\# \widehat W_\epsilon )}$ on $(\widehat W_\epsilon \setminus \cup _{p'\in\pi^{-1}(p)} B_{p'}(b_0))\subset \widehat M\#\widehat W_\epsilon$. By multiplying $g^{T(M\# W)}$ with a constant, we can and we will assume that $\delta=n^2$. 

Moreover,  we also see from the compactness of  $ \overline{ W  \setminus   B_{p}(b_0)}$ that (\ref{1.1}) still holds for $g^{T(\widehat M\# \widehat W_\epsilon )}$, possibly with a different $c$. 

By taking $\epsilon>0$ small enough, we see from (\ref{1.1}) and (\ref{1.2})  that $\widehat f:( \widehat M\#\widehat W_\epsilon, g^{T(\widehat M\# \widehat W_\epsilon )})\rightarrow S^n(1)$ verifies the   area decreasing condition, is locally constant near infinity, and that 
\begin{align}\label{1.3}
 k^{g^{T(\widehat M\# \widehat W_\epsilon )}}\geq n^2 \ \ {\rm on\ the\ support\ of} \  {\rm d}f.
\end{align} 
 By \cite[Theorems 2.1 and 2.2]{Z20}, one sees that if ${\rm deg}(\widehat f)={\rm deg}(f)\neq 0$, then 
${\rm inf}(k^{g^{T(\widehat M\# \widehat W_\epsilon )}})<0 $, which contradicts  the assumption that $k^{g^{T(M\#W)}}>0$.

$\ $

\noindent{\bf Acknowledgments.} The authors would like to thank Yuguang Shi and Guofang Wang for helpful discussions. W. Zhang was partially supported by NSFC Grant no. 11931007 and the Nankai Zhide Foundation. X. Wang was partially supported by NSFC Grant no. 12101361, the project of Young Scholars of SDU and the Fundamental Research Funds of SDU, Grant no. 2020GN063.

\def\cprime{$'$} \def\cprime{$'$}
\providecommand{\bysame}{\leavevmode\hbox to3em{\hrulefill}\thinspace}
\providecommand{\MR}{\relax\ifhmode\unskip\space\fi MR }
\providecommand{\MRhref}[2]{\href{http://www.ams.org/mathscinet-getitem?mr=#1}{#2}
}
\providecommand{\href}[2]{#2}

\end{document}